\documentstyle[amscd,amssymb,verbatim]{amsart}
\newcommand{\lrar}[1]{\begin{picture}(50,10)(-25,-5)
\put(-25,0){\vector(1,0){50}}
\put(0,5){\makebox(0,0)[b]{\mbox{$#1$}}}
\end{picture}}

\newcommand{\ldar}[1]{\begin{picture}(10,50)(-5,-25)
\put(0,25){\vector(0,-1){50}}
\put(5,0){\mbox{$#1$}}
\end{picture}}

\newcommand{\dil}{\operatorname{dil}}
\newcommand{\GL}{\operatorname{GL}}

\newcommand{\inv}{\operatorname{inv}}

\newcommand{\grad}{\operatorname{grad}}

\newcommand{\Tr}{\operatorname{Tr}}

\newcommand{\KK}{{\cal K}}
\newcommand{\RR}{{\cal R}}

\newcommand{\G}{{\Bbb G}}
\newcommand{\A}{{\Bbb A}}

\newcommand{\lan}{\langle}
\newcommand{\ran}{\rangle}

\newcommand{\Frob}{\operatorname{Frob}}

\renewcommand{\P}{{\Bbb P}}

\newcommand{\si}{\sigma}

\newcommand{\eps}{\epsilon}

\numberwithin{equation}{section}

\newtheorem{thm}{Theorem}[section]
\newtheorem{prop}[thm]{Proposition}
\newtheorem{lem}[thm]{Lemma}
\newtheorem{cor}[thm]{Corollary}
\newenvironment{rem}{\vspace{3mm}\noindent
{\bf Remark.}}{\vspace{3mm}}
\newenvironment{defi}{\vspace{3mm}\noindent
{\bf Definition.}}{\vspace{3mm}}
\newenvironment{rems}{\vspace{3mm}
\noindent {\bf Remarks.}}{\vspace{3mm}}

\newcommand{\Pf}{\noindent {\it Proof}}
\newcommand{\id}{\operatorname{id}}

\newcommand{\ov}{\overline}

\newcommand{\rk}{\operatorname{rk}}
\newcommand{\ra}{\rightarrow}

\newcommand{\FF}{{\cal F}}

\newcommand{\VV}{{\cal V}}

\newcommand{\LL}{{\cal L}}

\newcommand{\Hom}{\operatorname{Hom}}
\newcommand{\Ext}{\operatorname{Ext}}

\renewcommand{\a}{\alpha}
\renewcommand{\b}{\beta}
\newcommand{\om}{\omega}

\newcommand{\la}{\lambda}

\newcommand{\C}{{\Bbb C}}

\newcommand{\Z}{{\Bbb Z}}
\newcommand{\Q}{{\Bbb Q}}

\newcommand{\Ga}{\Gamma}

\newcommand{\wt}{\widetilde}

\newcommand{\sub}{\subset}
\newcommand{\ed}{\qed\vspace{3mm}}
\renewcommand{\check}{\vee}

\title{Generalized character sums associated to regular prehomogeneous
vector spaces}
\author{David Kazhdan and Alexander Polishchuk}
\thanks{Both authors are partially supported by NSF grants} 

\begin{document}
\maketitle

\bigskip

The purpose of this note is to give a short derivation of the
finite field analogue of Sato's functional equation for the zeta function
associated with a prehomogeneous vector space (see \cite{S}). 
We restrict ourselves to the case of a regular prehomogeneous
vector space,
however, we allow to twist our character sums by local systems
associated to arbitrary representations
of the component group of the stabilizer of a generic point.
The main idea of our approach is to use the Picard-Lefschetz
formula in $l$-adic cohomology instead of
using a lift of a prehomogeneous space to the characteristic zero
(as it is done in \cite{DG}).
Also we deduce another functional equation associated with
a regular prehomogeneous vector space (theorem \ref{thm3}).  

\section{Preliminaries and formulation of the main results}

Throughout this paper 
$k$ denotes an algebraically closed field of characteristic $p$,
$k_0$ denotes its finite subfield.
By a sheaf on a scheme $X$ over $k$ or $k_0$ we mean an object of
the derived category of constructible $l$-adic sheaves, where $l\neq p$  
(see e.g. \cite{D1}). For a morphism $f:X\ra Y$ we denote
by $f_!$, $f_*$, $f^!$, $f^*$ the natural derived functors between
the sheaves on $X$ and $Y$. Often we use the notation $F(f)$ for $f^*F$
where $F$ is a sheaf on $Y$.
For a sheaf $F$ we denote by
$\sideset{^p}{^i}{H}F$ the $i$-th perverse cohomology of $F$
(see \cite{BBD}).
For a scheme $X$ defined over $k_0$ we denote by $\Frob_{k_0}$
the geometric Frobenius morphism on $X$.
Recall that according to theorem of Laumon (see \cite{L0})
if $F$ is a sheaf on a $k$-scheme $X$ then one has the
equality of two kinds of Euler-Poincar\'e characteristics
$$\chi(X,F)=\chi_c(X,F).$$

By a multiplicative character of $k$ we always mean a continuous
homomorphism from $\hat{\Z}(1)(k)=\projlim_{(N,p)=1}\mu_N(k)$ 
to $\ov{\Q}_l^*$.
For every multiplicative
character $\chi$ we denote by $L_{\chi}$ the corresponding
Kummer sheaf on $\G_m$. 
We will often use the following theorem of Deligne (see \cite{I}):
for every sheaf $K$ on $\G_m$ and every multiplicative character $\chi$
one has
$\chi(\G_m,L_{\chi}\otimes K)=\chi(\G_m,K)$. 

If $\chi$ is a character of $k_0^*$ then it can be considered as a
multiplicative character of $k$ via
the natural homomorphism $\hat{\Z}(1)(k)\ra k_0^*$. Then 
the sheaf $L_{\chi}$ is defined over $k_0$ and $\Frob_{k_0}$
acts on the fiber of $L_{\chi}$ at $a\in\G_m(k_0)$ as the scalar
$\chi(a)$.

If $G$ is a connected algebraic group acting on $X$
then one can define the category of $G$-equivariant
perverse sheaves on $X$ (see e.g. \cite{Lu}). 
More generally, one can consider {\it relative $G$-equivariant} perverse
sheaves on $X$ corresponding to the pair $(\a,\chi)$ where
$\a:G\ra\G_m$ is a group homomorphism, $\chi$ is a
multiplicative character of $k$. These 
are perverse sheaves $F$ such that there exists an isomorphism
$$m^*F[\dim G]\simeq L_{\chi}(\a)\boxtimes F[\dim G]$$
of perverse sheaves on $G\times X$, 
where $m:G\times X\ra X$ is the action morphism (such 
an isomorphism on $G\times X$ can be chosen canonically by
imposing the condition that its restriction to $e\times X$ is
the identity, where $e\in G$ is the neutral element).
If the action of $G$ on $X$ is transitive,
and $F$ is a relative $G$-equivariant perverse sheaf on $X$ then
there exists a local system $\LL$ on $X$ such that 
$F=\LL[\dim X]$. One can show that a subquotient of
a relative $G$-equivariant perverse sheaf is again relative
$G$-equivariant (the proof is similar to \cite{Lu2},(1.9.1)).
In particular, an irreducible object in the category of
relative $G$-equivariant perverse sheaves is irreducible as a
perverse sheaf.

Let $\psi$ be a non-trivial additive character of
$k_0$. We denote by $L_{\psi}$ the corresponding
Artin-Schreier sheaf on $\A^1$ (defined over $k_0$).
Let $V$ be a finite-dimensional vector space over $k$. Then the
Fourier-Deligne transform (see \cite{L}) of a sheaf $F$ on $V$ is the
sheaf
$$\FF(F)=(p_2)_!(p_1^*F\otimes L_{\psi}(\lan v,v^{\vee}\ran))[\dim V]$$
on the dual space $V^{\vee}$, where $p_i$ are the projections
of $V\times V^{\vee}$ on two factors, $\lan v,v^{\vee}\ran$ denotes
the pairing between $V$ and $V^{\vee}$. 

For a multiplicative character $\chi$ we denote
$$G(\chi)=G(\chi,\psi)=H^1_c(\G_m,L_{\chi}\otimes L_{\psi}).$$
This is a 1-dimensional $\ov{\Q}_l$-space. If $\chi$ is a character
of $k_0^*$ then the action of the Frobenius on $G(\chi)$
is given by  
$$g(\chi)=g(\chi,\psi)=-\sum_{a\in k_0^*}\chi(a)\psi(a).$$ 
One has the following isomorphism of sheaves on $\G_m$:
\begin{equation}\label{isom2}
j_1^*\FF(j_{1!}L_{\chi})\simeq G(\chi)\otimes L_{\chi^{-1}}
\end{equation}
where $j_1:\G_m\ra\A^1$ is the standard
embedding. Furthermore, if $\chi$ is non-trivial then
\begin{equation}\label{isom2.5}
\FF(j_{1!}L_{\chi})\simeq G(\chi)\otimes j_{1!}L_{\chi^{-1}}.
\end{equation}

The facts about prehomogeneous spaces we use below are essentially
contained in \cite{KS}, we only need to make slight changes
due to finite characteristics (see e.g. section 1 of \cite{KKF}).    
Let $G$ be a connected reductive group over $k$,
$V$ be a {\it regular prehomogeneous vector space} over $k$ for $G$.
This means that $G$ acts linearly on $V$
with an open orbit $U\subset V$ and there is a relative invariant
$f\in k[V]$ with the character $\a\in\Hom(G,\G_m)$: $f(gx)=\a(g)f(x)$,
such that $F=\grad \log f$ gives a dominant map
from $U$ to $V^{\check}$. We fix such an invariant $f$ once and for all.
Often we consider $f$ as a morphism $U\ra\G_m$, in particular, we
consider sheaves of the form $L_{\chi}(f)$ on $U$ where $\chi$ is a 
multiplicative character.
The dual space $V^{\check}$ is also a regular prehomogeneous
space for $G$, $F$ induces a $G$-equivariant morphism of $U$ onto an open
$G$-orbit $U^{\check}\subset V^{\vee}$.
We assume that our prehomogeneous vector space is {\it separable} in the sense
of \cite{KKF},
i.e. the map $F$ is separable and the open orbit morphism 
$G\ra V^*:g\mapsto gv^{\vee}$ is separable for generic $v^{\vee}$.
Then as shown in \cite{KKF}
there exists a relative invariant $f^{\check}\in k[V^{\check}]$ with the
character $\a^{-1}$ and the map $F^{\check}=\grad\log f^{\check}$
is inverse to $F$.
Moreover, in this case the closed set $V\setminus U$
(resp. $V^{\vee}\setminus U^{\vee}$) coincides with the hypersurface $f=0$
(resp. $f^{\vee}=0$).
Since the function $f^{\vee}(F(x))f(x)$ on $U$
is $G$-invariant it should be constant, hence, we can rescale
$f^{\vee}$ in such a way that
\begin{equation}\label{eq1}
f^{\vee}(F(x))=f(x)^{-1}.
\end{equation}
Substituting $x=F^{\vee}(x^{\vee})$ into this equation we obtain
\begin{equation}\label{eq1.5}
f(F^{\vee}(x^{\vee}))=f^{\vee}(x)^{-1}.
\end{equation}
Since $F$ is homogeneous of degree $-1$ the equation (\ref{eq1}) also implies
that $\deg f^{\vee}=\deg f$.

Let $H$ be the stabilizer of a point $x_0\in U$.
Then $G$-equivariant local systems on $U$ correspond to representations
of the group $\ov{H}$ of connected components of $H$.
For an irreducible representation
$\rho$ of $\ov{H}$ let us denote by $\VV_{\rho}$ the corresponding
local system on $U$.

For a multiplicative character $\chi$ 
and an irreducible $\ov{H}$-representation $\rho$ we consider the 
relative $G$-equivariant perverse sheaf
$$S_V(\chi,\rho)=j_!(L_{\chi}(f)\otimes\VV_{\rho})[n]$$
on $V$ where $j:U\ra V$ is the natural open embedding,
$n=\dim V$.

We want to determine the Fourier-Deligne transform 
of the perverse sheaf 
$S_V(\chi,\rho)$ when it is irreducible.
We'll show below that this is the case for almost all $\chi$.
Also we'll show that the restriction of $\FF(S_V(\chi,\rho))$
to $U^{\check}$ is non-zero, therefore, by irreducibility
$\FF(S_V(\chi,\rho))$ is the Goresky-MacPherson extension of its
restriction to $U^{\check}$ which should be of the form
$L_{\chi^{-1}}(f^{\check})\otimes\VV^{\check}[n]$
for some irreducible $G$-equivariant sheaf $\VV^{\check}=\VV_{\rho^{\vee}}$
on $U^{\check}$.
It follows that for generic $\chi$ one has
$$\FF(S_V(\chi,\rho))\simeq S_{V^{\vee}}(\chi^{-1},\rho^{\check}).$$
Our main problem will be to determine $\rho^{\check}$. We'll show that
it doesn't depend on $\chi$ but only on $\rho$. Furthermore, notice
that there is a natural isomorphism between the stabilizer $H$ of a point
$x\in U$ and the stabilizer of the point $F(x)\in U^{\check}$
(since $F$ is $G$-equivariant). Thus, $\VV^{\check}$ corresponds
to some irreducible representation of $\ov{H}$ depending only on the original
representation $\rho$.
The explicit form of this dependence it given in the following theorem.

\begin{thm}\label{mainthm} 
There exists a homomorphism $\eps:\ov{H}\ra\{\pm 1\}$ such that \\
\noindent
(a) for almost all multiplicative characters $\chi$ one has 
\begin{equation}\label{isom0}
\FF(S_V(\chi,\rho))\simeq S_{V^{\vee}}(\chi^{-1},\eps\otimes\rho).
\end{equation}

\noindent
(b) for all characters $\chi$ one has
\begin{equation}\label{isom1}
(j^{\vee})^*\FF(S_V(\chi,\rho))\simeq (j^{\vee})^*
S_{V^{\vee}}(\chi^{-1},\eps\otimes\rho).
\end{equation}
\end{thm}

Now assume that our prehomogeneous vector space and a relative invariant $f$
are defined over a finite field $k_0$ and assume that
$U(k_0)$ is non-empty (then $f^{\vee}$ is also defined over $k_0$).
Then for every character $\chi$ of $k_0^*$ the sheaves $S_V(\chi,\rho)$ can be
equipped with the natural action of Frobenius but
the isomorphism (\ref{isom1}) is not in general compatible
with the action of Frobenius on LHS and RHS. The difference between two
actions is given by a non-zero constant depending on $\chi$.
We will show that this constant is essentially a product of some Gauss sums.
Then considering the traces of Frobenius acting on both sides of 
(\ref{isom1}) one derives some identities for character sums. 
Here is a more precise formulation.
Let us fix a point $x_0\in U(k_0)$. Then there is a natural
action of $\Frob_{k_0}$ on the component group $\ov{H}$. 
Let $\rho$ be an irreducible representation of $\ov{H}$ in a $\ov{\Q}_l$-space
$W$. 

\begin{defi} An action of $\Frob_{k_0}$ on the representation $\rho$ is
an operator $A\in\GL(W)$ such that 
$$\rho(\Frob_{k_0}(h))=A^{-1}\rho(h)A$$
and $A^N=\id$ for some $N>0$ 
\end{defi}

\begin{rem} An action of $\Frob_{k_0}$ on $\rho$ exists
if and only if $\ov{H}$-representations $\rho\circ\Frob_{k_0}$ and $\rho$ 
are isomorphic.
\end{rem}

Let us fix an action of $\Frob_{k_0}$ on $\rho$.
Then $A$ descends to an isomorphism
of $G$-equivariant local systems $\Frob_{k_0}^*\VV_{\rho}\ra\VV_{\rho}$. 
Notice that our assumption $A^N=\id$ implies that the Weil sheaf
$\VV_{\rho}$ is pure of weight $0$.
For every finite extension $k_1$ of $k_0$ and
$x\in U(k_1)$ we define a conjugacy class $\rho_x$ in
$\GL(W)$ as follows: choose an element $g\in G/H_0(\ov{k_0})$
such that $gx_0=x$ (where $H_0$ is the identity component of $H$),
then $\Frob_{k_1}(g)=gh$ for some $h\in H/H_0=\ov{H}$, now $\rho_x$
is the conjugacy class of $A^{[k_1:k_0]}\rho(h^{-1})$. Then we have
$$\Tr(\Frob_{k_1},(\VV_{\rho})_x)=\Tr(\rho_x).$$

We apply similar constructions to the dual prehomogeneous space $V^{\vee}$
taking as a fixed point $x_0^{\vee}=F(x_0)\in U^{\vee}(k_0)$.
Also we fix a square root of $|k_0|$ in $\ov{\Q}_l^*$. This leads to
a choice of a square root of $|k_1|$ for every finite extension 
$k_0\subset k_1$. Recall that an algebraic number $a\in\ov{\Q}_l^*$ is called 
pure
of weight $w$ relative to $|k_0|$ if for any isomorphism $i:\ov{\Q}_l\ra\C$
one has $|i(a)|=|k_0|^{w/2}$.

\begin{thm}\label{thm2}
For every irreducible representation $\rho$
of $\ov{H}$ equipped with an action of $\Frob_{k_0}$
there exists a finite extension $k'_0$ of $k_0$, two collections
of characters of $(k'_0)^*$: $(\la_1,\ldots,\la_{m+d})$ and
$(\mu_1,\ldots,\mu_m)$, where
$d=\deg(f)$, and a constant $\zeta\in\ov{\Q}_l^*$ of weight $0$, 
such that 
for every finite extension $k_1$ of $k'_0$ and
for every character $\chi$ of $k_1^*$ one has
\begin{eqnarray}\label{eq2}
(-1)^nq^{-\frac{n}{2}}
\sum_{x\in U(k_1)}\Tr(\rho_x)\chi(f(x))\psi(T\lan x,x^{\vee}\ran)=
\zeta^{[k_1:k_0]}\times
\nonumber\\
\prod_{i=1}^{m+d}\frac{g(\chi(\la_i\circ N))}{\sqrt{q}}\cdot
\prod_{j=1}^m\frac{g(\chi^{-1}(\mu_j\circ N))}{\sqrt{q}}\cdot
\Tr((\rho\otimes\eps)_{x^{\vee}})\chi^{-1}((-1)^mf^{\vee}(x^{\vee}))
\end{eqnarray}
for every $x^{\vee}\in U^{\vee}(k_1)$, where $q=|k_1|$,
$T:k_1\ra k_0$ is the trace, $N:k_1^*\ra (k'_0)^*$ is the norm.
Furthermore, if $\chi$ is different from all the
characters $\la_i^{-1}\circ N$ and $\mu_j\circ N$
then the sum in the LHS vanishes
for $x^{\vee}$ in the complement to $U^{\vee}$.
\end{thm}

\begin{rems} 
1. If all the characters $\la_i$ and $\mu_j$ are compositions of the norm
homomorphism and some characters of a subextension $k''_0\sub k'_0$ 
then one can replace $k'_0$ by $k''_0$ in the above formulation. 
  
2. In the case when $\rho$ is trivial and the prehomogeneous
space is obtained by reduction from characteristic zero,
according to \cite{DG} one has $m=0$ in the identity (\ref{eq2}).
Also in this case $\zeta$ is a root of unity.
We don't know whether this remains true for arbitrary representation
$\rho$. 
\end{rems}

\begin{cor} Let $\chi$ be one of the characters
$\la_i^{-1}$ or $\mu_j$. Then
$$\FF(j_{!*}(L_{\chi}(f)\otimes\VV_{\rho}[n]))|_{U^{\check}}=0.$$
\end{cor}

\Pf . This follows easily from the fact that the weight of the RHS in 
(\ref{eq2}) drops for such $\chi$. Indeed, in this case
the pure perverse sheaf of weight $2n$
$$\FF(j_{!*}(L_{\chi}(f)\otimes\VV_{\rho}[n]))|_{U^{\check}}$$
is a quotient of the perverse sheaf
$$\FF(j_!(L_{\chi}(f)\otimes\VV_{\rho}[n]))|_{U^{\check}}$$ 
which is pure of weight $<2n$, hence the former sheaf should be zero.
\ed

As before we fix a square root of $|k_0|$ in $\ov{\Q}_l^*$. This allows  
to define
the square root of the Tate twist $F\mapsto F(\frac{1}{2})$ which multiplies
the action of $\Frob_{k_0}$ by $|k_0|^{-\frac{1}{2}}$.

\begin{thm}\label{thm3}
Keeping the notations of theorem \ref{thm2} for every character $\chi$
of $(k'_0)^*$ one has the following 
isomorphism of sheaves over $k'_0$:
\begin{align*}
&\FF \wt{j}_{!*}
(L_{\psi}(\frac{f(x)u_1\ldots u_m}{t_1\ldots t_{m+d}})\otimes
L_{\chi}(\frac{f(x)u_1\ldots u_m}{t_1\ldots t_{m+d}})\otimes
\VV_{\rho}(x)\otimes\\
&\bigotimes_i L_{\la_i}(t_i^{-1})\otimes
\bigotimes_j L_{\mu_j}(u_j^{-1}))(m+\frac{d+n}{2})
\simeq\zeta^{\deg}\otimes\wt{j}^{\vee}_{!*}
(L_{\psi}(\frac{(-1)^{m+d} t^{\vee}_1\ldots t^{\vee}_{m+d}}
{f^{\vee}(x^{\vee})u^{\vee}_1\ldots u^{\vee}_m})
\otimes\\
&L_{\chi}(\frac{(-1)^{m+d} t^{\vee}_1\ldots t^{\vee}_{m+d}}
{f^{\vee}(x^{\vee})u^{\vee}_1\ldots u^{\vee}_m})\otimes
\VV_{\rho\eps}(x^{\vee})\otimes
\bigotimes_i L_{\la_i}(-t^{\vee}_i)\otimes\bigotimes_j 
L_{\mu_j}(-u^{\vee}_j)) 
\end{align*}
where $\wt{j}:U\times\G_m^{2m+d}\ra V\times\A^{2m+d}$ (resp.
$\wt{j}^{\vee}:U^{\vee}\times\G_m^{2m+d}\ra V^{\vee}\times\A^{2m+d}$) is the
natural open embedding, 
$(x,t_1,\ldots,t_{m+d},u_1,\ldots,u_m)$ are the coordinates on 
$V\times\A^{2m+d}$,
$(x^{\vee},t^{\vee}_1,\ldots,t^{\vee}_{m+d},u^{\vee}_1,\ldots,u^{\vee}_m)$ 
are the dual coordinates on $V^{\vee}\times\A^{2m+d}$,
$\zeta^{\deg}$ is the geometrically constant sheaf of rank $1$ on which
$\Frob_{k_1}$ acts as $\zeta^{[k_1:k_0]}$,
$\wt{j}_{!*}$ denotes the Goresky-MacPherson extension (conjugated by
an appropriate shift).
\end{thm}

\begin{rem} In the situation of the above theorem one can construct also the 
following isomorphism of sheaves:
\begin{align*}
&\FF \wt{j}_{!*}
(L_{\psi}(\frac{f(x)u_1\ldots u_mt_{m+d}}{t_1\ldots t_{m+d-1}})\otimes
L_{\chi}(\frac{f(x)u_1\ldots u_m}{t_1\ldots t_{m+d-1}})\otimes\\
&\VV_{\rho}(x)\otimes\bigotimes_{i=1}^{m+d-1} L_{\la_i}(t_i^{-1})\otimes
\bigotimes_{j=1}^m L_{\mu_j}(u_j^{-1}))(m+\frac{d+n}{2})
\simeq\zeta^{\deg}\otimes\\
&\wt{j}^{\vee}_{!*}
(L_{\psi}(\frac{(-1)^{m+d} f^{\vee}(x^{\vee})u^{\vee}_1\ldots u^{\vee}_m
t^{\vee}_{m+d}}{t^{\vee}_1\ldots t^{\vee}_{m+d-1}})\otimes
L_{\la_{m+d}}(\frac{(-1)^{m+d} f^{\vee}(x^{\vee})u^{\vee}_1\ldots u^{\vee}_m
t^{\vee}_{m+d}}{t^{\vee}_1\ldots t^{\vee}_{m+d-1}})\otimes\\
&\VV_{\rho\eps}(x^{\vee})\otimes
\bigotimes_{i=1}^{m+d-1} L_{\la_i}(-t^{\vee}_i)\otimes\bigotimes_{j=1}^m 
L_{\mu_j}(-u^{\vee}_j) \otimes L_{\chi}(-t^{\vee}_{m+d}))
\end{align*}
The proof is completely analogous to that of theorem \ref{thm3}.

\end{rem}

\section{Auxiliary results}

\begin{lem}\label{lem2} For almost all $\chi$ the 
perverse sheaf $S_V(\chi,\rho)$ is irreducible.
\end{lem}

\Pf . Consider the closed embedding 
$i:U\ra V\times\G_m$ such that $i(x)=(x,f(x))$. Let 
$K=i_*\VV_{\rho}$. Then according
to \cite{GL} Theorem 2.3.1 (based on \cite{KL} (6.5.2)) the canonical map
$$p_{1!}(K\otimes p_2^*L_{\chi})\ra p_{1*}(K\otimes p_2^*L_{\chi})$$
is an isomorphism for almost all $\chi$. Hence, the canonical map
$$j_!(L_{\chi}(f)\otimes\VV_{\rho})[n]\ra
j_*(L_{\chi}(f)\otimes\VV_{\rho})[n]$$
is an isomorphism for almost all $\chi$. Therefore, for almost all $\chi$
the sheaf $S_V(\chi,\rho)$ coincides with the Goresky-MacPherson 
extension of an irreducible perverse sheaf
$L_{\chi}(f)\otimes\VV_{\rho}[n]$.
\ed

\begin{rem} To find explicitly the finite set of multiplicative characters
$\chi$ for which $S_V(\chi,\rho)$ is not irreducible is a difficult problem.
When $\rho=1$ this set can be found from the Bernstein polynomial of $f$
(see \cite{DG}). In general one can at least assert that
if $S_V(\chi,\rho)$ is not irreducible then there exists a point 
$x\in V\setminus U$
such that $(\a|_{H_x})^*L_{\chi}$ is trivial, where $H_x$ is the stabilizer 
of $x$.
In the cases when $G$-orbits on $V$ are known this allows to give a bound
on the set of exceptional characters.
\end{rem}

\begin{lem}\label{lem1} 
Let $\VV$ be an irreducible $G$-equivariant local system on $U$. 
Then there exists a multiplicative character of finite
order $\chi_{\VV}$ and a canonical
isomorphism of $G$-equivariant sheaves on $\G_m\times U$ (where
$G$ acts on the second factor):
$$\dil^*\VV\simeq L_{\chi_{\VV}}\boxtimes\VV$$
where $\dil:\G_m\times U\ra U$ maps $(t,x)$ to $tx$.
\end{lem}

\Pf. The local system $\dil^*\VV$ on $\G_m\times U$ is $G$-equivariant
with respect to the action of $G$ on the second factor, hence
we have 
$$\dil^*\VV\simeq L\boxtimes\VV'$$ 
for some local system $L$ on 
$\G_m$ and some $G$-equivariant local system $\VV'$ on $U$. Restricting
this isomorphism to $\{1\}\times U$ we find that $\VV'\simeq\VV$ and
the rank of $L$ is equal to $1$.
It remains to prove that 
$L$ has form $L_{\chi_{\VV}}$. To this end let
us fix a point $x_0\in U$ and let $H_{l}$ be the stabilizer of the line
$l=l_{x_0}$ spanned by $x_0$.
Then we have a surjective homomorphism $\pi:H_{l}\ra\G_m$, 
hence there exists a homomorphism $\si:\G_m\ra H_{l}$
such that $\pi\circ\si:\G_m\ra\G_m$ is of the form $t\mapsto t^N$
for some $N>0$.  
Now from $G$-equivariance we derive that
the pull-back of $L$ by $\pi$ is trivial. 
Hence, $[\pi\circ\si]^*L$ is trivial, so $L\simeq L_{\chi_{\VV}}$ 
for some $\chi_{\VV}$ of finite order.
\ed

\begin{lem}\label{lem2.5} Let $K$ be a sheaf on $\G_m$ such that
$H^i(\G_m,L_{\chi}\otimes K)=0$ for every multiplicative character $\chi$
and every $i\neq 0$. Then $K$ is perverse.
\end{lem}

\Pf .
Since the functor
$R\Ga_c(\G_m,?)$ has cohomological amplitude $[0,1]$ with respect
to the perverse $t$-structure our assumption on $K$ immediately implies
that $\sideset{^p}{^i}{H}K=0$ for $i\neq -1, 0$, while for
any $\chi$ one has
$$H^0_c(\G_m,L_{\chi}\otimes\sideset{^p}{^{-1}}{H}K)=
H^1_c(\G_m,L_{\chi}\otimes\sideset{^p}{^0}{H}K)=0$$
and
$$H^0_c(\G_m,L_{\chi}\otimes K)=
H^1_c(\G_m,L_{\chi}\otimes\sideset{^p}{^{-1}}{H}K)\oplus
H^0_c(\G_m,L_{\chi}\otimes\sideset{^p}{^0}{H}K)$$
This implies that $\chi(\G_m,\sideset{^p}{^{-1}}{H}K)=0$,
hence $R\Ga_c(\G_m,L_{\chi}\otimes\sideset{^p}{^{-1}}{H}K)=0$
for every $\chi$. Applying theorem of Laumon (see
Proposition 3.4.5 of \cite{GL}) we obtain that 
$\sideset{^p}{^{-1}}{H}K=0$, so $K$ is perverse.
\ed

We say that a sheaf $K$ on $\G_m$ is $!$-hypergeometric
if it is obtained as a multiplicative shift of the
$!$-convolution of a finite number
of sheaves of the form $L_{\chi}\otimes L_{\psi}|_{\G_m}[1]$
or $\inv^*(L_{\chi}\otimes L_{\psi}|_{\G_m})[1]$ where
$\inv:\G_m\ra\G_m:t\mapsto t^{-1}$. 
Recall (see e.g. Proposition 8.5.2 of \cite{K})
that any perverse sheaf on $\G_m$ has non-negative
Euler-Poincar\'e characteristic and the only simple perverse
sheaves with zero Euler-Poincar\'e characteristic are the sheaves of the form
$L_{\chi}[1]$. By a theorem of Katz every simple perverse
sheaf with Euler-Poincar\'e characteristic equal to $1$ is $!$-hypergeometric.
The following proposition is a slight strengthening of this result.

\begin{prop}\label{prop}
Let $K$ be a sheaf on $\G_m$ such that
$\chi(\G_m,K)=1$ and $H^i_c(\G_m,L_{\chi}\otimes K)=0$ for $i\neq 0$
for every multiplicative character $\chi$.
Then $K$ is a $!$-hypergeometric sheaf.
\end{prop}

\Pf . It follows from lemma \ref{lem2.5} that $K$ is perverse. 
We use induction in length of $K$. If $K$ is simple then
the assertion follows from Theorem 8.5.3 of \cite{K}.
Assume that $K$ is not simple and let $K_0$ be a simple perverse
subobject of $K$, so that we have an exact sequence
\begin{equation}\label{ex}
0\ra K_0\ra K\ra K_1\ra 0.
\end{equation}
We necessarily have $\chi(\G_m,K_1)\le 1$.
Suppose $\chi(\G_m,K_1)=0$. Then $K_1$ is a successive extension of sheaves
of the form $L_{\chi}[1]$. However, for any $\chi$ one has
$$\Hom(K,L_{\chi}[1])\simeq H^0(\G_m,D(K\otimes L_{\chi^{-1}}[1]))\simeq
H^0_c(\G_m,K\otimes L_{\chi^{-1}}[1])^*=0$$
(where $D$ denotes the Verdier duality) which is a contradiction. 
Therefore, $\chi(\G_m,K_1)=1$. Furthermore,
for any $\chi$ the space
$\Hom(K_1,L_{\chi}[1])$ is a subspace of $\Hom(K,L_{\chi}[1])=0$. Thus,
$H^1_c(\G_m, K_1\otimes L_{\chi^{-1}})=0$ and by induction assumption
$K_1$ is a $!$-hypergeometric sheaf. On the other hand,
$\chi(\G_m,K_0)=0$, hence $K_0\simeq L_{\chi_0}[1]$ for some character 
$\chi_0$.
Notice that the exact sequence (\ref{ex}) doesn't split (otherwise,
we would have $H^1_c(\G_m,L_{\chi_0^{-1}}\otimes K)\neq 0$). It remains
to use the fact that every non-trivial extension of a $!$-hypergeometric
sheaf $K_1$ by $L_{\chi_0}[1]$ is $!$-hypergeometric which follows
from the fact that $\dim\Ext^1(K_1,L_{\chi_0}[1])=1$ and 
Theorem 8.4.7 of \cite{K}.
\ed

\begin{lem}\label{lem8} Assume that for some
collections of characters $\la_i$, $\mu_j$, $\la'_k$, $\mu'_l$
defined over finite field $k_0$ one has the identity
$$\prod g(\chi(\la_i\circ N))g(\chi^{-1}(\mu_j\circ N))=
c^{[k_1:k_0]}\chi(a)\prod g(\chi(\la'_k\circ N))g(\chi^{-1}(\mu'_l\circ N))$$
for every multiplicative character $\chi$ of every
finite extension $k_1$ of $k_0$, where $c\in\ov{\Q}_l^*$,
$a\in k_0^*$ are some constants, $N:k_1^*\ra k_0^*$ is the norm.
Then the collections $(\la_i)$ and $(\la'_k)$ (resp.
$(\mu_j)$ and $(\mu'_l)$) coincide.
\end{lem}

\Pf . Both parts of the identity are multiplicative Fourier transform
of trace functions of $!$-hypergeometric sheaves on $\G_m$.
Hence, by Theorem 1.1.2 of \cite{L} semisimplifications of these sheaves
are isomorphic. Looking at the monodromy at $0$ and $\infty$ we conclude
that the corresponding collections of characters should coincide
(see \cite{K}, Theorems 8.4.11 and 8.4.12).
\ed

\begin{lem} \label{lem7} Let $\chi$ be a non-trivial multiplicative
character. Then for every sheaf $K$ on $\A^1$ one has a 
canonical isomorphism
$$H^0_c(\G_m,L_{\chi}\otimes j_1^*\FF(K))\simeq
G(\chi)\otimes H^0_c(\G_m,L_{\chi^{-1}}\otimes j_1^*K).$$
\end{lem}

\Pf . One has
$$H^0_c(\G_m,L_{\chi}\otimes j_1^*\FF(K))\simeq
H^0_c(\G_m\times\A^1, L_{\chi}(a)\otimes L_{\psi}(at)\otimes
p_2^*K[1])$$
where $a$ is the coordinate on $\G_m$, $t$ is the coordinate on $\A^1$.
Projecting to the second factor and using (\ref{isom2.5}) we obtain
the result.
\ed

\section{Computation}

To compute the Fourier transform of $j_!(L_{\chi}(f)\otimes V_{\rho})$
we want to combine the information coming from $G$-equivariance and
from the stationary phase approximation. Technically the latter is
obtained from the Picard-Lefschetz formula. 
The idea is to use the condition that
$F=\grad\log f$ induces a birational isomorphism from $V$ to $V^{\vee}$.
Let us define a rational function
on $V\times\A^{d}$ where $d=\deg f$ by
$$R(x,t_1,\ldots,t_{d})=\frac{f(x)}{t_1\ldots t_{d}}.$$
We will denote by $y=(x,t_1,\ldots,t_{d})$ the coordinate in
$V\times\A^{d}$ and by $\xi=(x^{\vee},t_1^{\vee},\ldots,t_{d}^{\vee})$ 
the coordinate in the dual space $V^{\vee}\times\A^{d}$. We denote
$$\lan y,\xi\ran=\lan x,x^{\vee}\ran + t_1t_1^{\vee}+\ldots+
t_{d}t_{d}^{\vee}.$$

\begin{lem}\label{lem3} The map $y\mapsto \grad R(y)$ induces an isomorphism
from $U\times\G_m^{d}$ to $U^{\vee}\times\G_m^{d}$ with the inverse
$\xi\mapsto (-1)^d R^{\vee}(\xi)^{-2}\grad R^{\vee}(\xi)$
where $R^{\vee}$ is defined in a similar
way as $R$ starting from $f^{\vee}$ instead of $f$.
\end{lem}

\Pf . This follows easily from (\ref{eq1}). 
\ed

The above lemma implies that for $\xi\in U^{\vee}\times\G_m^{d}$
the function
\begin{equation}\label{Rxi}
R_{\xi}(y)=R(y)+\lan y,\xi\ran
\end{equation}
has a unique critical
point in $U\times\G_m^{d}$, namely $y_{\xi}=(\grad R)^{-1}(\xi)$. 
Moreover, $y_{\xi}$
is a non-degenerate quadratic singularity of the corresponding
level set of $R_{\xi}$ (see \cite{R}, sec.~6; in particular,
if characteristics of $k$ is equal to $2$ this implies that
$n+d$ is even).

\begin{lem}\label{lem3.5} 
The unique critical value of $R_{\xi}$ is equal to
$(-1)^dR^{\check}(\xi)^{-1}$.
\end{lem}

\Pf . According to lemma \ref{lem3} the unique critical point of 
$R_{\xi}$ is
$$y_{\xi}=(-1)^dR^{\vee}(\xi)^{-2}\grad R^{\vee}(\xi),$$
hence
$$\lan y_{\xi}, \xi\ran=(-1)^dR^{\vee}(\xi)^{-2}
\cdot\lan \grad R^{\vee}(\xi),\xi\ran=0$$
since $\deg R^{\vee}=0$. Hence,
the critical value is equal to
$$R_{\xi}(y_{\xi})=R(y_{\xi})=R(\grad R^{\vee}(\xi)).$$
Now using (\ref{eq1.5}) it is easy to see that
$$R(\grad R^{\vee}(\xi))=(-1)^d R^{\vee}(\xi)^{-1}.$$
\ed

\begin{lem}\label{lem4} For every multiplicative character $\chi$,
an irreducible $G$-equivariant local system $\VV$ on $U$ and a point
$\xi\in U^{\vee}\times\G_m^{d}$ one has
$$G(\chi^{-1})\otimes\FF(J_!(L_{\chi}(R)\otimes p_U^*\VV))|_{\xi}\simeq
R\Ga_c(\G_m,L_{\chi^{-1}\chi_{\VV}}\otimes
j_1^*\FF(R_{\xi!}(p_U^*\VV))[n+d])$$
where $J:U\times\G_m^{d}\ra V\times\A^{d}$
is the natural embedding,
$R_{\xi}:U\times\G_m^{d}\ra\A^1$ is the morphism defined by (\ref{Rxi}),
$p_U$ is the projection to $U$, the character
$\chi_{\VV}$ was defined in lemma \ref{lem1}, $j_1:\G_m\ra\A^1$ is the
standard embedding.
\end{lem}
      
\Pf .
Using the isomorphism (\ref{isom2}) we can write
\begin{align*}
&G(\chi^{-1})\otimes\FF(J_!(L_{\chi}(R)\otimes p_U^*\VV))|_{\xi}\simeq\\
&\simeq R\Ga_c(\G_m\times U\times\G_m^{d},
L_{\chi^{-1}}(a)\otimes L_{\psi}(aR(y)+\lan y,\xi\ran)\otimes
p_U^*\VV[n+d+1])
\end{align*}
where $a$ is coordinate on the first factor $\G_m$.
Now making a change of variables $y\mapsto ay$ and using
the fact that $R$ has degree $0$ and lemma \ref{lem1} we obtain
\begin{align*}
&R\Ga_c(\G_m\times U\times\G_m^{d},
L_{\chi^{-1}}(a)\otimes L_{\psi}(aR(y)+\lan y,\xi\ran)\otimes
p_U^*\VV[n+d+1])\\
&\simeq
R\Ga_c(\G_m\times U\times\G_m^{d},
L_{\chi^{-1}\chi_{\VV}}(a)\otimes L_{\psi}(aR_{\xi}(y))
\otimes p_U^*\VV)[n+d+1])\\
&\simeq
R\Ga_c(\G_m\times\A^1, p_1^*L_{\chi^{-1}\chi_{\VV}}\otimes L_{\psi}(at)\otimes
p_2^*R_{\xi!}(p_U^*\VV)[n+d+1])\\
&\simeq
R\Ga_c(\G_m,L_{\chi^{-1}\chi_{\VV}}\otimes j_1^*\FF(R_{\xi!}(p_U^*\VV))[n+d]).
\end{align*}
\ed

\begin{lem}\label{lem5} For every $\xi\in U^{\vee}\times\G_m^{d}$
the sheaf $j_1^*\FF(R_{\xi!}(p_U^*\VV))[n+d]$ on $\G_m$ is perverse.
\end{lem}

\Pf . The sheaf $(J^{\vee})^*\FF(J_!(L_{\chi}(R)\otimes p_U^*\VV))[n+d]$ on
$U^{\vee}\times\G_m^d$ is perverse and relative equivariant with respect
to the transitive action of $G\times\G_m$. Therefore,
it is equal to some local system shifted by $[n+d]$. Now 
lemma \ref{lem4} implies 
that for any $\chi$ the cohomology of $R\Ga_{c}(\G_m,L_{\chi}\otimes
j_1^*\FF(R_{\xi!}(p_U^*\VV))[n+d])$ is concentrated in degree zero.
It remains to apply lemma \ref{lem2.5}.
\ed

\begin{lem} \label{lem6} For every $\xi\in U^{\vee}\times\G_m^{d}$ and
every $i\neq n+d$ the sheaf $\sideset{^p}{^i}{H} R_{\xi!}(p_U^*\VV)$
is constant (up to a shift).
\end{lem}

\Pf . This follows from lemma \ref{lem5} and from $t$-exactness
of the Fourier transform (see \cite{L}).
\ed

Let $\VV=\VV_{\rho}$ be
an irreducible $G$-equivariant local system on $U$,
$\chi$ be a multiplicative character. 
The sheaf $F=j_!(L_{\chi}(f)\otimes\VV)$ on $V$ is relative
$G$-equivariant with characters $(\a,\chi)$. Therefore,
its Fourier transform is also relative $G$-equivariant with the
same characters. Taking into account the fact that $\FF(F)[n]$ is perverse
we conclude that
$$(j^{\vee})^*\FF(j_!(L_{\chi}(f)\otimes\VV))
\simeq L_{\chi^{-1}}(f^{\vee})\otimes \VV^{\check}$$
where $\VV^{\check}$ is a $G$-equivariant sheaf on $U^{\vee}$.
Therefore, we have
$$(J^{\vee})^*\FF(J_!(L_{\chi}(R)\otimes p_U^*\VV))\simeq
G(\chi^{-1})^{d}\otimes
L_{\chi^{-1}}(R^{\vee})\otimes p_{U^{\vee}}^*\VV^{\vee}$$
where $p_{U^{\vee}}:U^{\vee}\times\G_m^d\ra U^{\vee}$ is the projection.
Let us fix a point $\xi\in U^{\vee}\times\G_m^{d}$ such that
$R^{\vee}(\xi)=1$.
Then we have an isomorphism of $G(\chi^{-1})^{d}\otimes
\VV^{\vee}_{\xi}$, where $\VV^{\vee}_{\xi}$ is the fiber of 
$p_{U^{\vee}}^*\VV^{\vee}$ at $\xi$,
with the fiber of $\FF(J_!(L_{\chi}(R)\otimes p_U^*\VV))$ at $\xi$.
Using lemma \ref{lem4} we can write
\begin{equation}\label{isom3}
\VV^{\vee}_{\xi}\simeq
G(\chi^{-1})^{-d-1}
\otimes R\Ga_c(\G_m, L_{\chi^{-1}\chi_{\VV}}\otimes
j_1^*\FF(R_{\xi!}(p_U^*\VV))[n+d]).
\end{equation}
Notice that this isomorphism is compatible with the grading 
(where $\deg G(\chi^{-1})=0$) so the cohomology of the RHS is
concentrated in degree $0$. 

\vspace{3mm}
\noindent 
{\it Proof of theorem \ref{mainthm}}. 
Assume that $\chi\neq\chi_{\VV}$. Then using lemma
\ref{lem7} and lemma \ref{lem6}
we deduce from (\ref{isom3}) a (non-canonical) isomorphism
\begin{equation}\label{isom4}
\VV^{\vee}_{\xi}\simeq H^0_c(\G_m, L_{\chi\chi^{-1}_{\VV}}\otimes K).
\end{equation}
where $K=j_1^*\sideset{^p}{^{n+d}}{H}R_{\xi!}(p_U^*\VV)$.
Let $H_{\xi}\subset G\times\G_m^d$ be the stabilizer of $\xi$,
$\ov{H}_{\xi}$ be its group of connected components.
The isomorphism (\ref{isom4}) is compatible with the action of $\ov{H}_{\xi}$
(the action on the RHS is induced by the 
action of $H_{\xi}$ on $V\times\G_m^{d}$ which preserves the function
$R_{\xi}$).
We can write $K=\oplus_{W} W\otimes K_W$
where the sum is taken over irreducible representations of $\ov{H}_{\xi}$,
$K_W$ are perverse sheaves on $\G_m$. 
Then we have
$$\VV^{\vee}_{\xi}\simeq\oplus_W W^{\oplus\chi_c(\G_m,K_W)},$$
According to Grothendieck-Ogg-Shafarevich formula we have
$$\chi_c(\G_m,K_W)=-\sum_{t\in\G_m(k)} a_t(K_W)-s_0(K_W)-
s_{\infty}(K_W)$$
where for every $t\in \P^1(k)$ we denote by $s_t(F)$ the Swan conductor
of $F$ at $t$, $a_t(F)=\rk(F)+s_t(F)-\rk(F|_t)$.
Note that all the numbers $-a_t(K_W)$ are non-negative since $K_W$ are
perverse (see lemma (2.2.1.1) of \cite{L}). 
Furthermore, we have
$$
-a_t(K_W)\ge\rk(K_W[-1])-\rk(K_W|_t[-1])=\rk\Phi_t(K_W[-1])
$$
where $\Phi_t$ denotes the vanishing cycles functor near the point $t$
(indeed, the restriction of $K_W[-1]$ to some non-empty open subset
of $\G_m$ is a local system concentrated in degree $0$, 
hence $s_t(K_W[-1])\ge 0$). 
Thus, we have an isomorphism of $\ov{H}_{\xi}$-representations
$$\VV^{\vee}_{\xi}\simeq \oplus_{t\in\G_m} \Phi_t(K[-1])\oplus P$$
for some $\ov{H}_{\xi}$-representation $P$.
Note that by lemma \ref{lem6} we have
$$\Phi_t(K[-1])=\Phi_t(R_{\xi!}(p_U^*\VV)[n+d-1])$$
According to lemmas \ref{lem3} and \ref{lem3.5} 
the function $R_{\xi}$ has a unique critical
point $y_0$ and the corresponding critical value is $t_0=(-1)^d$. 
Since $y_0$ is a unique critical point of $R_{\xi}$
on $U\times\G_m^{d}$, it is preserved by the action of $H_{\xi}$. 
Hence the $H_{\xi}$-module $\Phi_{t_0}(R_{\xi!}(p_U^*\VV)[n+d-1])$ 
contains
$$\Phi(p_U^*\VV[n'-1])|_{y_0}
\simeq \VV_{y_0}\otimes \Phi(\ov{\Q}_l[n+d-1])|_{y_0}$$ 
as a direct summand, where $\Phi$ denotes
the functor of vanishing cycles with respect to $R_{\xi}$,
$\VV_{y_0}$ is the fiber of $p_U^*\VV$ at $y_0$.
Furthermore, since $y_0$ is a non-degenerate singular point
we have according to Picard-Lefschetz formula (see \cite{D})
that $L=\Phi(\ov{\Q}_l[n+d-1])|_{y_0}$ is a one-dimensional
vector space concentrated in degree $0$.
Note that the isomorphism of $H_{\xi}$ with $H_{y_0}$ induced by the
action of $H_{\xi}$ on $U\times\G_m^{d}$ coincides with the
isomorphism between $H_{p_{U^{\vee}}(\xi)}$ and $H_{p_U(y_0)}$
considered before. Let us denote by $H$ any of these group.
The $H$-representation
$\VV^{\vee}_{\xi}$ contains $\VV_{y_0}\otimes L$ as a direct summand.
The action of $H$ on $\VV_{y_0}$ is given by the representation
$\rho$ where $\VV=\VV_{\rho}$. On the other hand,
the space of vanishing cycles $L$ has a canonical generator
up to a sign. This implies that the action of $H$ on it is given by
some character $\eps:H\ra\{\pm 1\}$.
According to lemma 2 for almost all $\chi$ the local system $\VV^{\vee}$
is irreducible, hence $\VV^{\vee}_{\xi}$ coincides with its direct
summand $\VV_{y_0}\otimes L$, so we have 
$\VV^{\vee}\simeq \VV_{\rho\otimes\eps}$.
This proves part (a) of Theorem (\ref{mainthm}).

Now invoking the isomorphism (\ref{isom3}) which holds
for all characters $\chi$ we deduce that
\begin{equation}\label{M}
j_1^*\FF(R_{\xi!}(p_U^*\VV))[n+d]=(\rho\otimes\eps)\otimes M
\end{equation}
for some sheaf $M$ on $\G_m$ such that
$R\Ga_c(\G_m,L_{\chi}\otimes M)$ is
concentrated in degree zero for all $\chi$, and $\chi(\G_m,M)=1$.
By lemma \ref{lem5} $M$ is perverse, hence,
applying proposition \ref{prop} we conclude that
$M$ is a $!$-hypergeometric sheaf.
Thus, for every character $\chi$ 
we have a canonical isomorphism
\begin{equation}\label{isom5}
(j^{\vee})^*\FF j_!(L_{\chi}(f)\otimes\VV_{\rho})\simeq
G(\chi^{-1})^{-d-1}\otimes H^0_c(\G_m, L_{\chi^{-1}\chi_{\VV}}\otimes M)
\otimes L_{\chi^{-1}}(f^{\vee})\otimes \VV_{\rho\otimes\eps}.
\end{equation}
This finishes the proof of theorem (\ref{mainthm}).

An additional information that we get from the above argument is 
recorded in the following lemma. 

\begin{lem}\label{lem9}
The sheaf $R_{\xi!}(p_U^*\VV)$ is tamely ramified on $\P^1$ and
is smooth on $\G_m\setminus\{(-1)^d\}$.
\end{lem}

\vspace{3mm}
\noindent
{\it Proof of theorem \ref{thm2}}.
Let $M$ be the $!$-hypergeometric sheaf defined by (\ref{M}). If $\xi$ is
a $k_0$-point then $M$ is defined over $k_0$. It follows that
for some $a\in k^*$ the sheaf $[a]^*M$ (where $[a](b)=ab$) is isomorphic
(forgetting the action of Frobenius) to 
the $!$-convolution of sheaves of the form $L_{\la_i}\otimes L_{\psi}[1]$,
$i=1,\ldots,m_1$, $\inv^*(L_{\mu_j}\otimes L_{\psi})[1]$,
$j=1,\ldots,m_2$, where $\la_i$ and $\mu_j$ are multiplicative
characters of some finite extension $k'_0$ of $k_0$. 
First we claim that $m_1=m_2+1$ and $a=(-1)^{m_2+d}$. Indeed,
lemma \ref{lem9} implies that
$\FF(j_{1!}M)$ is tamely ramified on $\P^1$ and is smooth outside
$(-1)^{d+1}$.
But $\FF(j_{1!}M)|_{\G_m}$ is the $!$-hypergeometric sheaf
obtained by a multiplicative shift by $a^{-1}$ of the $!$-convolution of
the sheaves $L_{\psi}[1]$, $L_{\mu_j}\otimes L_{\psi}[1]$, $j=1,\ldots m_2$,
and $\inv^*(L_{\la_i}\otimes L_{\psi})[1]$, $i=1,\ldots, m_1$.
According to Theorem 8.4.11 of \cite{K} such a sheaf is tamely ramified
only if $m_1=m_2+1$, and smooth outside $(-1)^{d+1}$ only if $(-1)^{m_1}a=
(-1)^{d+1}$, which proves our claim. 
Taking traces of Frobenius acting on both sides
of (\ref{isom5}) we get the equality of the form
\begin{eqnarray}\label{eq3}
\FF_{k_1}j_!(\Tr(\rho_x)\chi(f(x)))(x^{\vee})=
c^{[k_1:k_0]}\cdot\chi(-1)^{m_2+d}\times
\nonumber\\
\frac{\prod_{i=1}^{m_2+1}g(\chi^{-1}(\la'_i\circ N))
\prod_{j=1}^{m_2}
g(\chi(\mu'_j\circ N))}{g(\chi^{-1})^{d+1}}\cdot
\Tr((\rho\otimes\eps)_{x^{\vee}})\chi^{-1}(f^{\vee}(x^{\vee}))
\end{eqnarray}
for every $x^{\vee}\in U^{\vee}(k_1)$ where $k_1$ is a finite extension
of $k'_0$, $N:k_1^*\ra (k'_0)^*$ is the norm, 
$\FF_{k_1}$ is the discrete Fourier transform on $V(k_1)$
defined as follows:
\begin{equation}\label{Four}
\FF_{k_1}(\phi)(x^{\vee})=(-1)^nq^{-\frac{n}{2}}
\sum_{x\in V(k_1)}\phi(x)\psi(T\lan x,x^{\vee}\ran)
\end{equation}
where $q=|k_1|$, $T:k_1\ra k_0$ is the trace. 
Note that for generic $\chi$ the sheaf
$j_!L_{\chi}(f)\otimes\VV_{\rho}$ is pure of weight $0$.
Since the Fourier-Deligne transform sends pure sheaves to pure sheaves
raising weights by $n$ we find from (\ref{eq3})
that the constant $c$ should be of weight $d-2m_2$ relative to $|k_0|$.

It remains to get rid of the denominator in the RHS of (\ref{eq3}), i.e.
to prove that $d+1$ of the characters $\la'_i$ are trivial.
The idea is to compute the Fourier transform by a different method
to obtain an equality similar to (\ref{eq3}) but with the denominator in the 
RHS only containing Gauss sums of the form $g(\chi^{-1}\chi_i)$ with
$\chi_i$ non-trivial. Then lemma \ref{lem8} will imply that
the denominator in (\ref{eq3}) can be cancelled.

Let us fix non-trivial multiplicative characters $\chi_1,\ldots,\chi_d$
of $k'_0$ and consider the smooth sheaf
$$\wt{\VV}=\VV\boxtimes  L_{\chi_1}\boxtimes\ldots L_{\chi_d}$$
on $U\times\G_m^d$.
We want to compute the Fourier transform of the sheaf 
$J_!(L_{\chi}(R)\otimes\wt{\VV})$ on $V_d=V\times\A^d$ using the Radon
transform (see e.g. Lemma 3.3.9 of \cite{DG}). Since 
$R$ has degree zero we have
$$\FF J_!(L_{\chi}(R)\otimes\wt{\VV})\simeq 
\RR J_!(L_{\chi}(R)\otimes\wt{\VV}):=p_{2!}(p_1^* J_!
(L_{\chi}(R)\otimes\wt{\VV})\otimes\om)[n+d]$$
where $\om$ is the sheaf on $V_d\times V_d^{\vee}$ defined as
the cone of the natural adjunction morphism 
$e_!e^!\ov{\Q}_l\ra\ov{\Q}_l$ for
the closed embedding $e$ of the affine
hyperplane $Q=\{(y,\xi):\lan y,\xi\ran=1\}$ into 
$V_d\times V_d^{\vee}$.
Note that the stabilizer $H_{\xi}$ of $\xi\in V_d^{\vee}$ acts on
$V_d$ preserving $Q\cap p_2^{-1}(\xi)$, hence the sheaf
$\om_{\xi}:=\om|_{p_2^{-1}(\xi)}$ is
$H_{\xi}$-equivariant and the isomorphism
$$\FF(J_!(L_{\chi}(R)\otimes\wt{\VV}))|_{\xi}\simeq 
R\Ga_c(U\times\G_m^d,L_{\chi}(R)\otimes\wt{\VV}\otimes J^*\om_{\xi}[n+d])$$
is compatible with  the $H_{\xi}$-action.
We can write
$$R\Ga_c(U\times\G_m^d, L_{\chi}(R)\otimes\wt{\VV}\otimes 
J^*\om_{\xi}[n+d])\simeq
R\Ga_c(\G_m,L_{\chi}(R)\otimes\KK_{\xi})$$
where $\KK_{\xi}=R_!(\wt{\VV}\otimes J^*\om_{\xi}[n+d])$.
Choosing $\xi$ in the form $(x^{\vee};1,\ldots,1)$ where
$f^{\check}(x^{\vee})=1$ we obtain that
$$\FF(j_!(L_{\chi}(f)\otimes\VV))|_{x^{\vee}}\simeq
\bigotimes_{i=1}^d G(\chi^{-1}\chi_i)^{-1}\otimes R\Ga_c(\G_m,\KK_{\xi}).$$
We have an $H_{\xi}$-isomorphism 
$\KK_{\xi}\simeq (\rho\otimes\eps)\otimes K_{\xi}$
and by proposition \ref{prop} the sheaf
$K_{\xi}$ is $!$-hypergeometric.
Thus, we obtain the equality similar to (\ref{eq3}) but with constant
factor of the form
$$(\wt{c})^{[k_1:k_0]}\cdot\chi(\wt{a})\cdot
\frac{\prod_{i=1}^{\wt{m}_1}g(\chi^{-1}(\wt{\la}_i\circ N))
\prod_{j=1}^{\wt{m}_2}g(\chi(\wt{\mu}_j\circ N))}{\prod_{i=1}^d 
g(\chi^{-1}(\chi_i\circ N))}$$
Since the characters $\chi_i$ are non-trivial we can finish
the proof of (\ref{eq2}) by applying lemma \ref{lem8}.

It remains to prove the last statement of Theorem \ref{thm2}.
Let us fix $\rho$. For every finite extension $k'_0\sub k_1$
and a character $\chi$ of $k_1^*$ we can write the equality (\ref{eq2}) in
the form
$$(j^{\vee})^*\FF_{k_1}(j_!\phi_{\chi})=C_{\chi}\cdot\phi^{\vee}_{\chi}$$
where 
$$\phi_{\chi}(x)=\Tr(\rho_x)\chi(f(x))$$ 
for $x\in U(k_1)$,
$$\phi^{\vee}_{\chi}(x^{\vee})=\Tr((\rho\otimes\eps)_{x^{\vee}})\chi^{-1}((-1)^m
f^{\vee}(x^{\vee}))$$ 
for $x^{\vee}\in U^{\vee}(k_1)$, 
$$C_{\chi}=\zeta^{[k_1:k_0]}\cdot
\prod_{i=1}^{m+d}\frac{g(\chi(\la_i\circ N))}{\sqrt{q}}\cdot
\prod_{j=1}^m\frac{g(\chi^{-1}(\mu_j\circ N))}{\sqrt{q}}.$$
We want to prove that if $\chi$ is different from all the characters
$\la_i\circ N$, $\mu_j\circ N$ then $\FF_{k_1}(j_!\phi_{\chi})$
is zero outside $U^{\vee}(k_1)$. First we claim that one can
assume $k_1$ to be sufficiently large. Indeed, assume that
there exists a finite extension $k_1\sub k'_1$ such that for every
finite extension $k'_1\subset k_2$ the function 
$\FF_{k_2}(j_!\phi_{\chi\circ N_{k_2/k_1}})$ vanishes outside 
$U^{\vee}(k_2)$. This implies that the sheaf $\FF(j_!L_{\chi})$
is zero outside $U^{\vee}$ (by Theorem 1.1.2 of \cite{L}), 
in particular, $\FF_{k_1}(j_!\phi_{\chi})$ vanishes
outside $U^{\vee}(k_1)$. 
Now we 
can choose $k_1$ large enough so that there exists
a character $\chi_0$ of $k_1^*$ (different from all the characters
$\la_i\circ N$ and $\mu_j\circ N$) such that 
$$\FF_{k_1}(j_!\phi_{\chi_0})=C_{\chi_0}\cdot j^{\vee}_!\phi^{\vee}_{\chi_0}$$
(this follows from theorem \ref{mainthm} (a)).
Using the identity $||\FF_{k_1}(\phi)||=||\phi||$ we derive that
\begin{equation}\label{chig}
||j_!\phi_{\chi_0}||=|C_{\chi_0}|\cdot ||j^{\vee}_!\phi^{\vee}_{\chi_0}||
\end{equation}
Now the statement follows imeediately from the fact 
that the norms $||j_!\phi_{\chi}||$ and $||j^{\vee}_!\phi^{\vee}_{\chi}||$ 
do not depend on $\chi$. 
Indeed, assume that $\chi$ is different from all the characters $\la_i\circ N$
and $\mu_j\circ N$. Then we have $|C_{\chi}|=|C_{\chi_0}|$, so
we can replace $\chi_0$ by $\chi$ in the equality (\ref{chig}).
It follows that the norm of the function
$\FF_{k_1}(j_!\phi_{\chi})$ is equal to the norm of its restriction to
$U^{\vee}(k_1)$, hence this function 
vanishes outside $U^{\vee}(k_1)$.
\ed

\vspace{3mm}
\noindent
{\it Proof of theorem \ref{thm3}}.
By irreducibility it suffices to prove the isomorphism of restrictions
of sheaves in both sides to $U^{\vee}\times\G_m^{2m+d}$. 
The idea is to combine (\ref{eq2}) with the simple identity
\begin{equation}\label{eqsim}
(1-q)\cdot\psi(a)\la(a)=\sum_{\nu}g(\la\nu)\nu^{-1}(a)
\end{equation}
where $a\in k_1^*$, $\la$ is a character of $k_1^*$,
the sum is taken over all characters of $k_1^*$, by abuse of notation
we denote by $\psi$ the additive character of $k_1$ composed of the
trace and $\psi$.
Applying (\ref{eqsim}) we get
\begin{align*}
& (1-q)\cdot (\wt{j}^{\vee})^*
\FF_{k_1}\wt{j}_!(\psi(\frac{f(x)u_1\ldots u_m}{t_1\ldots t_{m+d}})
(\chi\circ N)(\frac{f(x)u_1\ldots u_m}{t_1\ldots t_{m+d}})
\Tr(\rho(r_x))\times\\
&\prod_i \la_i^{-1}(N t_i)\prod_j \mu_j^{-1}(N u_j))=
\sum_{\nu} g(\nu(\chi\circ N))\cdot
(\wt{j}^{\vee})^*\FF_{k_1}\wt{j}_!(\nu^{-1}(f(x))\Tr(\rho(r_x))\times\\
&\prod_i (\nu(\la_i^{-1}\circ N))(t_i)
\prod_j(\nu(\mu_j\circ N))^{-1}(u_j)).
\end{align*}
Here $\FF_{k_1}$ denotes the Fourier transform on $(V\times\A^{2m+d})(k_1)$
defined by a formula similar to (\ref{Four}) but with $n$ replaced
by $n+2m+d$, $N:k_1^*\ra (k'_0)^*$ is the norm.

Now applying (\ref{eq2}) we can rewrite this as
\begin{align*}
& \zeta^{[k_1:k_0]}\cdot q^{-2m-d}\cdot
\sum_{\nu}g(\nu(\chi\circ N))\nu((-1)^m\frac{f^{\vee}(x^{\vee})
u^{\vee}_1\ldots u^{\vee}_m}{t^{\vee}_1\ldots t^{\vee}_{m+d}})
\Tr((\rho\eps)(r_{x^{\vee}}))\times \\
&\prod_i\la_i (N t^{\vee}_i)g(\nu^{-1}(\la_i\circ N))
g(\nu(\la_i^{-1}\circ N))
\prod_j\mu_j (N u^{\vee}_j)g(\nu(\mu_j\circ N))
g(\nu^{-1}(\mu_j\circ N)^{-1})=\\
&\zeta^{[k_1:k_0]}\cdot\sum_{\nu}g(\nu(\chi\circ N))
\nu(\frac{(-1)^{m+d} f^{\vee}(x^{\vee})
u^{\vee}_1\ldots u^{\vee}_m}{t^{\vee}_1\ldots t^{\vee}_{m+d}})
\Tr((\rho\eps)(r_{x^{\vee}}))\times\\
&\prod_i\la_i(N(-t^{\vee}_i))\prod_j\mu_j(N(-u^{\vee}_j))+
(1-q)\cdot\sum_{i=1}^M m_i t_{k_1}(F_i)
\end{align*}
where the coefficients $m_i\in\Z$ do not depend on the extension
$k_1$, $t_{k_1}(F_i)$ are the trace functions of some simple
perverse sheaves  
which are relative equivariant with respect to the natural action of
$\G_m\times\G_m^{2m+d}$ on $U\times\G_m^{2m+d}$ (where the
first factor $\G_m$ acts on $U$ by dilations). Finally we can use
(\ref{eqsim}) again to deduce the following equality
\begin{align*}
&(\wt{j}^{\vee})^*\FF_{k_1}
\wt{j}_!(\psi(\frac{f(x)u_1\ldots u_m}{t_1\ldots t_{m+d}})
(\chi\circ N)(\frac{f(x)u_1\ldots u_m}{t_1\ldots t_{m+d}})
\Tr(\rho(r_x))\times\\
&\prod_i \la_i^{-1}(N t_i)\prod_j\mu_j^{-1}(Nu_j))=\\
&=\zeta^{[k_1:k_0]}\psi(\frac{(-1)^{m+d}t^{\vee}_1\ldots t^{\vee}_{m+d}}
{f^{\vee}(x^{\vee})
u^{\vee}_1\ldots u^{\vee}_m})
(\chi\circ N)(\frac{(-1)^{m+d}t^{\vee}_1\ldots t^{\vee}_{m+d}}
{f^{\vee}(x^{\vee})u^{\vee}_1\ldots u^{\vee}_m})
\Tr((\rho\eps)(r_{x^{\vee}}))\times\\
&\prod_i\la_i(N(-t^{\vee}_i))\prod_j\mu_j(N(-u^{\vee}_j))+
\sum_{i=1}^M m_it_{k_1}(F_i).
\end{align*}
The LHS of this identity is the trace function of the sheaf
$$(\wt{j}^{\vee})^*\FF\wt{j}_!
(L_{\psi}(\frac{f(x)u_1\ldots u_m}{t_1\ldots t_{m+d}})\otimes\LL)$$
where 
$$\LL=L_{\chi}(\frac{f(x)u_1\ldots u_m}{t_1\ldots t_{m+d}})\otimes
\VV_{\rho}(x)\otimes
\bigotimes_i L_{\la_i^{-1}}(t_i)\otimes\bigotimes_j L_{\mu_j^{-1}}(u_j).$$
Hence, it differs from the trace function of the sheaf
$$F=(\wt{j}^{\vee})^*\FF\wt{j}_{!*}
(L_{\psi}(\frac{f(x)u_1\ldots u_m}{t_1\ldots t_{m+d}})\otimes\LL)$$
by a linear combination of trace functions of simple perverse sheaves
(up to shift) on $U^{\vee}\times\G^{2m+d}$
which are not isomorphic to $F$. Assume that $F\neq 0$.
Then it is a simple perverse sheaf (up to shift), hence 
by Theorem 1.1.2 of \cite{L} it  coincides with one
of the sheaves corresponding to terms in the RHS.
It remains to prove that $F\neq 0$ and $F$ is not relative equivariant
with respect to the action of $\G_m\times\G_m^{2m+d}$.
Note that the sheaf $F$ is the image of the morphism of perverse sheaves
$F_!\ra F_*$ where $F_!$ (resp. $F_*$) is defined by the same formula
but using the $!$-extension (resp. $*$-extension).
Let us consider the action of $\G_m$ on $V\times\A^{2m+d}$
(resp. on $V^{\vee}\times\A^{2m+d}$) such that for 
$t\in\G_m$ 
the corresponding operator $\b(t):V\times\A^{2m+d}\ra V\times\A^{2m+d}$ 
(resp. $\b^{\vee}(t):V^{\vee}\times\A^{2m+d}\ra V^{\vee}\times\A^{2m+d}$) 
multiplies the coordinate $t_1$ (resp. $t^{\vee}_1$) by $t$ and leaves all the 
other coordinates unchanged. 
Let us denote by $Q$ the map from $U\times\G_m^{2m+d}$ to $\G_m$ given by
$$Q(x,t_1,\ldots,t_{m+d},u_1,\ldots,u_m)=
\frac{f(x)u_1\ldots u_m}{t_1\ldots t_{m+d}}.$$
Let us fix $\xi\in U^{\vee}\times\G_m^{2m+d}$ and denote
by $Q_{\xi}:U\times \G_m^{2m+d}\ra\A^1$ the map given by
$$Q_{\xi}(\wt{x})=Q(\wt{x})+\lan \wt{x},\xi\ran.$$
Then we claim that there exist a commutative diagram 
\begin{equation}\label{comm}
\begin{array}{ccc}
F_!(t\b^{\vee}(t)\xi)\otimes L_{\la_1^{-1}}&\lrar{} 
&\FF((Q_{\xi})_!\LL)|_{\G_m}\\
\ldar{}&&\ldar{}\\
F_*(t\b^{\vee}(t)\xi)\otimes L_{\la_1^{-1}}&\lrar{}
&\FF((Q_{\xi})_*\LL)|_{\G_m}
\end{array}
\end{equation}
where rows are isomorphisms, $t$ is the coordinate on $\G_m$.
Indeed, the first row is obtained
from the following sequence of isomorphisms:
\begin{align*}
&\FF((Q_{\xi})_!\LL)|_t\simeq 
H^*_c(U\times\G_m^{2m+d},L_{\psi}(tQ_{\xi}(\wt{x}))\otimes\LL(\wt{x}))\simeq\\
&\simeq
H^*_c(U\times\G_m^{2m+d},L_{\psi}(Q(\b(t)^{-1}\wt{x})+\lan \wt{x}, t\xi\ran)
\otimes\LL(\wt{x}))\simeq\\
&\simeq H^*_c(U\times\G_m^{2m+d},L_{\psi}(Q(\wt{x})+\lan \b(t)\wt{x}, t\xi\ran)
\otimes\LL(\b(t)\wt{x}))\simeq\\
&\simeq
H^*_c(U\times\G_m^{2m+d},L_{\psi}(Q(\wt{x})+\lan \wt{x}, t\b^{\vee}(t)\xi\ran)
\otimes\LL(\wt{x}))\otimes L_{\la_1^{-1}}|_t\simeq\\
&\simeq F_!(t\b^{\vee}(t)\xi)\otimes L_{\la_1^{-1}}|_t.
\end{align*}
Similarly, one constructs the lower arrow in the diagram (\ref{comm}) and
the commutativity can be easily checked.
Thus, it suffices to prove that the image of the morphism
$$(Q_{\xi})_!\LL\ra (Q_{\xi})_*\LL$$ is not
a smooth sheaf on $\G_m$. By exactness of the vanishing cycles
functor it suffices to prove that for some point $t_0\in\G_m$
the natural map 
\begin{equation}\label{map}
\Phi_{t_0}(Q_{\xi})_!\LL\ra\Phi_{t_0}(Q_{\xi})_*\LL
\end{equation}
is non-zero, where $\Phi_{t_0}$ denotes the vanishing cycles functor
at $t_0$. This follows from the fact that there exists
a unique non-degenerate critical point of $Q_{\xi}$ in
$U\times\G_m^{2m+d}$ with non-zero critical value $t_0$ (this can be checked
as in lemmas \ref{lem3} and \ref{lem3.5}) and Picard-Lefschetz formula. 
Indeed, according to this formula the non-degenerate critical point
gives a one-dimensional direct summand in both $\Phi_{t_0}((Q_{\xi})_!\LL)$
and $\Phi_{t_0}((Q_{\xi})_*\LL)$ which maps identically under the morphism
(\ref{map}).
\ed

\noindent
{\it Acknowledgment}. We are grateful to Pavel Etingof for helpful
discussions.


\begin{thebibliography}{99}
\bibitem{BBD} A.~Beilinson, J.~Bernstein, P.~Deligne, {\it Faisceaux
pervers}, Asterisque 100 (1982), 7--172.
\bibitem{D1} P. Deligne, {\it La conjecture de Weil, II},
Publ. Math. IHES 52 (1980), 313--428.
\bibitem{D} P. Deligne, {\it La formule de Picard-Lefschetz},
exp. XV in
{\it Groupes de monodromie en g\'eom\'etrie alg\'ebrique}
(SGA 7 II). Lecture Notes in Mathematics, Vol. 340.
Springer-Verlag, Berlin-New York, 1973.
\bibitem{DG} J. Denef, A. Gyoja, {\it Character sums
associated to prehomogeneous vector spaces}, Compositio Math. 113
(1998), 273--346.
\bibitem{GL} O. Gabber, F. Loeser, {\it Faisceaux pervers $l$-adique
sur un tore}, Duke Math. Journal 83 (1996), 501--606.
\bibitem{I} L. Illusie, {\it Th\'eorie de Brauer et charact\'eristique
d'Euler-Poincar\'e (d'apr\`es P. Deligne)}, in {\it
Caract\'eristique d'Euler-Poincar\'e}, Asterisque 82-83 (1981), 161--172.
\bibitem{K} N. Katz, {\it Exponential sums and differential equations},
Princeton University Press, 1990.
\bibitem{KL} N. Katz, G. Laumon, {\it Transformation de Fourier
et majoration de sommes exponentielles}, Publ. IHES 62 (1986), 361--418.
\bibitem{KKF} T. Kimura, T. Kogiso, M. Fujinaga,
{\it Fundamental theorem of prehomogeneous vector spaces of
characteristic $p$}, Bull. Austral. Math. Soc. 56 (1997), 331--341. 
\bibitem{KS} T. Kimura, M. Sato, {\it A classification of    
irreducible prehomogeneous vector spaces and their relative 
invariants}, Nagoya Math. J. 65 (1977), 1--155. 
\bibitem{L0} G. Laumon, {\it Comparaison de caract\'eristiques
d'Euler-Poincar\'e en cohomologie $l$-adique}, C.R. Acad. Sci. Paris
Ser. I Math. 292 (1981), 209--212.
\bibitem{L} G. Laumon, {\it Transformation de Fourier, constantes
d'\'equations fonctionnelles et conjecture de Weil},
Publ. IHES 65 (1987), 131--210.  
\bibitem{Lu} G.~Lusztig, {\it Intersection cohomology complexes on a
reductive group}, Invent. Math. 75 (1984), 205--272.
\bibitem{Lu2} G.~Lusztig, {\it Character sheaves. I}, Adv. in Math. 56 (1985),
193--237.
\bibitem{R} D. S. Rim, {\it Formal deformation theory}, exp. VI in
{\it Groupes de monodromie en g\'eom\'etrie alg\'ebrique}
(SGA 7 I). Lecture Notes in Mathematics, Vol. 288.
Springer-Verlag, Berlin-New York, 1972
\bibitem{S} M. Sato, {\it Theory of prehomogeneous vector spaces
(algebraic part)---the English translation of 
Sato's lecture from Shintani's note}.
Notes by T.~Shintani. Translated from the Japanese by M.~Muro.
Nagoya Math. J. 120 (1990), 1--34.
\end{thebibliography}
\end{document}